\documentclass{amsart} 

\usepackage{latexsym,amssymb,amscd,eucal}

\newtheorem{thm}{Theorem}[section]
\newtheorem{lem}[thm]{Lemma}
\newtheorem{prop}[thm]{Proposition}
\newtheorem{cor}[thm]{Corollary} 
\newtheorem{Def}[thm]{Definition}

\newcommand{\BQ}{{\mathbb{Q}}}
\newcommand{\BR}{{\mathbb{R}}}

\newcommand{\fC}{{\mathfrak{C}}}

\newcommand{\cc}{{\mathcal{C}}}
\newcommand{\Si}{{\Sigma}}
\newcommand{\ra}{{\rightarrow}}

\DeclareMathOperator{\Ann}{Ann}

\DeclareMathOperator{\id}{id}

\begin{document}

\title [the even cobordism category]
        {The parity of the maslov index and  the even cobordism category}
\author{ Patrick M. Gilmer}
\author{Khaled Qazaqzeh}
\address{Department of Mathematics\\
   Louisiana State University \\
   Baton Rouge, LA 70803-4918}
\email{gilmer@math.lsu.edu, kqazaq1@lsu.edu}
\urladdr{www.math.lsu.edu/\textasciitilde gilmer/}
\thanks{partially supported by NSF-DMS-0203486}

\keywords{even cobordisms, Lagrangian, Maslov index}

\date{April 12, 2005}

\begin{abstract}
 We give a formula for the parity of the Maslov index of a triple of Lagrangian subspaces of a skew symmetric bilinear form over $\BR$.  We define an index two subcategory (the even  subcategory)  of a 3-dimensional cobordism category. The objects of the category are surfaces are equipped with Lagrangian subspaces of their real first homology. This generalizes a result of the first author where surfaces are equipped with Lagrangian subspaces of their rational first homology. 
\end{abstract}

\maketitle
 
\section {introduction}
 In \cite{G}, the first author considered a cobordism category $\cc.$ This category can be described roughly as follows.   The objects of $\cc$ are closed surfaces equipped with Lagrangian subspaces of their rational first homology. A morphism of $\cc$ between $N:\Si \ra \Si'$ is a cobordism between  $ \Si$ and $\Si'$.   Also, the first author defined a subcategory $\cc^+$ of $\cc$ of index two. It would be more consistent with other work \cite{T} \cite{W} to consider a similarly defined cobordism category $\fC$ where the extra data of a Lagrangian subspace is a subspace of the real first homology. 
 The main goal of this article is define  an analogous index two subcategory $\fC^+$ of $\fC$.  
 We call $\fC^+$
   the even cobordism category.
    If one restricts
to this `index two' cobordism subcategory, one may
obtain functors, related to the TQFT functors defined by Turaev with initial data a modular category, but without taking a
quadratic extension of the ground ring of the modular category as
is sometimes needed in \cite[p.76]{T}.

 It is not possible to simply modify the proof given in \cite{G} for the existence of  $\cc^+$ to obtain a proof for the existence of $\fC^+$. This is because not every real Lagrangian subspace can  be realized as the kernel of the map induced on first homology by the inclusion of a surface to a 3-manifold which has the surface as its boundary. Only the subspaces which are completions of subspaces of the rational homology can be so realized. So another approach has to be used. We actually reduce the problem to the one already solved for $\cc$ but this requires some new algebraic results.  These algebraic results may be of independent interest. 
 
  We prove the algebraic results  in $\S2$. This section is written without any appeal to topology. It can be read independently of the rest of the paper.  We  derive the following congruence for the Maslov index, denoted $\mu$: 
  
  \begin{thm} \label{maincong}
Let $V$ be a symplectic vector space and $\lambda_{1},\lambda_{2}$, and
$\lambda_{3}$ be any three Lagrangian subspaces, then we have
\begin{align}\label{E:fivth} \notag
\mu(\lambda_{1},\lambda_{2},\lambda_{3}) &\equiv  \dim (\lambda_{1}) +
\sum_{1 \leq i < j \leq 3} \dim (\lambda_{i}\cap \lambda_{j})\quad \bmod(2)\\
\notag &\equiv \dim (\lambda_{1}) + \sum_{1 \leq i < j \leq 3}
\dim (\lambda_{i} + \lambda_{j}) \quad \bmod(2)
\end{align}
\end{thm}

 If  $\lambda_{1} \cap \lambda_{2}=\lambda_{2} \cap \lambda_{3}=\lambda_{1} \cap \lambda_{3}=0$, this result follows from  \cite[1.5.7]{LV} which gives a formula for  the Maslov index in terms of a special form these Lagrangians must take in this case.
We give a very diffferent proof.
 Theorem \ref{maincong}  will be the key to proving that the morphisms
of
$\fC^+$ are closed under composition. 
   In \S3, we describe the weighted cobordism categories $\cc$ and $\fC$ in greater detail.  In \S4, 
   $\fC^+$  is defined.  

\section{lagrangian subspaces and the maslov index }

 Let $V$ be a symplectic  vector space, i.e.  $V$ is  finite dimensional
over $\mathbb{R}$  and endowed with a skew symmetric bilinear form $\psi$. We do not require that the form is nondegenerate.  If $A$ is  a subspace of $V$, its annihilator,  $\Ann(A)$,  is  the set of elements which pair under the form with all of $A$ to give zero.  If $A$ and $A'$ are two subspaces, then \cite[IV.3.1.a, IV.3.1.1 ]{T} 
\begin{equation} \label{ann} \Ann (A + A') =  \Ann (A) \cap \Ann(A') \end{equation}
\begin{equation} \label{ann2} \Ann (A \cap A') =  \Ann (A) + \Ann(A') \end{equation}
A
subspace $A\subset V$ is said to be a Lagrangian subspace if $A$ =
$\Ann(A$).  We have not been able to find the following result in the literature.

\begin{thm}\label{T:firstt}
Let  $(V, \psi )$  be a symplectic vector space and
$\lambda_{1},\lambda_{2}$, and $\lambda_{3}$ be three Lagrangian subspaces.
Then we have
\begin{align} \notag
\dim (\lambda_{1} + \lambda_{2} + \lambda_{3}) \equiv
\dim (\lambda_{1}\cap\lambda_{2}\cap\lambda_{3})\quad \bmod(2)
\end{align}
\end{thm}

\begin{proof} We have an skew symmetric bilinear form $\psi$ on $V$. Now
define a form $\{\, , \}$ on
$(\lambda_{1}+\lambda_{2}+\lambda_{3})/(\lambda_{1}\cap\lambda_{2}\cap\lambda_{3})$
  by
$\{a,b\}=\psi(a,b)$ where $a,b\in
(\lambda_{1}+\lambda_{2}+\lambda_{3})/(\lambda_{1}\cap\lambda_{2}\cap\lambda_{3})$.
To show that this new form is well-defined, let $a_{1}, a_{2}\in
(\lambda_{1}+\lambda_{2}+\lambda_{3})$ such that $\bar{a_{1}}=\bar{a_{2}}$,
i.e.  $a_{1}-a_{2}\in (\lambda_{1}\cap\lambda_{2}\cap\lambda_{3})$. It
follows $\psi(a_{1}-a_{2},b)=0$ for all $b\in
(\lambda_{1}+\lambda_{2}+\lambda_{3})$, so $\psi(a_{1},b)=\psi(a_{2},b) $.
Hence $\{a_{1},b\}=\{a_{2},b\}$ for all $b\in
(\lambda_{1}+\lambda_{2}+\lambda_{3})$ that $\{\, , \}$ is well-defined. Since $\psi$ is skew symmetric
bilinear form, so is $\{\, , \}$.
We now wish to show that   $\{\, , \}$  is
non-degenerate. So let
$a\in(\lambda_{1}+\lambda_{2}+\lambda_{3})/(\lambda_{1}\cap\lambda_{2}\cap\lambda_{3})$
 such that $\{a,b\}=0$ for all
$b\in(\lambda_{1}+\lambda_{2}+\lambda_{3})/(\lambda_{1}\cap\lambda_{2}\cap\lambda_{3})$,
 i.e. $\psi(a,b)=0$ for all $b\in  \lambda_{1}+\lambda_{2}+\lambda_{3}$, it
implies that $a\in \Ann(\lambda_{1}+\lambda_{2}+\lambda_{3}$).  By equation \eqref{ann} 
\begin{align}
\notag
\Ann(\lambda_{1}+\lambda_{2}+\lambda_{3}) &= \Ann(\lambda_{1} +
\lambda_{2})\cap \Ann(\lambda_{3})\\ \notag & = (\Ann(\lambda_{1})\cap
\Ann(\lambda_{2})) \cap \lambda_{3}\\ \notag & =
\lambda_{1}\cap\lambda_{2}\cap\lambda_{3}.
\end{align}

So $a\in\lambda_{1}\cap\lambda_{2}\cap\lambda_{3}$, i.e. $a=0$ in
$(\lambda_{1}+\lambda_{2}+\lambda_{3})/(\lambda_{1}\cap\lambda_{2}\cap\lambda_{3})$.
As is well-known,  A non-degenerate symplectic vector space must be even dimensional. Hence
$(\lambda_{1}+\lambda_{2}+\lambda_{3})/(\lambda_{1}\cap\lambda_{2}\cap\lambda_{3})$
 is of even dimension, so we get
\begin{align}
\notag \dim (\lambda_{1} + \lambda_{2} + \lambda_{3}) \equiv
\dim (\lambda_{1}\cap\lambda_{2}\cap\lambda_{3}) \quad \bmod(2)
\end{align}
\end{proof}
 We have the following well-known proposition \cite[IV.3.5]{T}
\begin{prop}
Let $\lambda_{1}$, $\lambda_{2}$ and $\lambda_{3}$ be three Lagrangian
subspaces of $V$. Define a bilinear form $\langle\, , \rangle$ on
$(\lambda_{1}+\lambda_{2})\cap \lambda_{3}$ by
\begin{align} \label{E:second}
\langle a,b\rangle=\psi(a_{2},b)
\end{align} where $a, b\in (\lambda_{1}+\lambda_{2})\cap \lambda_{3}$ and
$a=a_{1}+a_{2}$.  $\langle\, , \rangle$ is a well-defined symmetric bilinear form.
\end{prop}

\begin{proof} To show is $\langle\, , \rangle$ is well-defined, note that the decomposition
$a=a_{1}+a_{2}$, where $a_1 \in A_1$ and  $a_2 \in A_2$,  is unique up to an element in $\lambda_{1}\cap\lambda_{2}$,
and this element annihilates $b$ for all $b\in\lambda_{1}+\lambda_{2}$. So
the form is well-defined.  As $\psi$ is bilinear, 
$\langle\, , \rangle$ is bilinear.

Let $a$ be as before and $b=b_{1}+b_{2}$ where
$b_{1}\in\lambda_{1}, b_{2}\in\lambda_{2}$ and $b\in\lambda_{3}$. Since
$\lambda_{i} = \Ann(\lambda_{i})$ for $i$ = 1, 2, 3  and $\psi$ is
skew symmetric, we have
\begin{align}
\notag \psi(a_{2},b)&=\psi(a-a_{1},b)& \\
\notag             &=\psi(a,b)-\psi(a_{1},b) \\
\notag             &=\psi(b,a_{1})  \\
\notag             &=\psi(b_{1}+b_{2},a_{1}) \\
\notag             &=\psi(b_{1},a_{1})+\psi(b_{2},a_{1})+\psi(b_{2},a_{2})\\
\notag             &=\psi(b_{2},a).
\end{align}
 Hence the form is symmetric.
\end{proof}

 \begin{Def}  The  Maslov index
$\mu(\lambda_{1},\lambda_{2},\lambda_{3})$ of the triple
$(\lambda_{1},\lambda_{2},\lambda_{3})$ is the signature of the form $\langle\, , \rangle$ defined
above. \end{Def}

In general, $\langle\, , \rangle$ is degenerate. In fact, it is known that its
annihilator contains  $(\lambda_{1}\cap\lambda_{3})$ +
$(\lambda_{2}\cap\lambda_{3})$ \cite[p.182-183]{T}. 
 If  $\lambda_1 \cap \lambda_2=0$, it is known that the annhilator  is $(\lambda_{1}\cap\lambda_{3})$ +
$(\lambda_{2}\cap\lambda_{3})$  \cite[1.5.6]{LV}. We show this is true in general.

\begin{thm}\label{T:secondt}
Let ($V$,$\psi$) be a symplectic vector space and $\lambda_{1}$,
$\lambda_{2}$, and $\lambda_{3}$ be three Lagrangian subspaces, then the
induced form $\langle\, , \rangle$ on ($\lambda_{1} + \lambda_{2}$) $\cap
\lambda_{3}$ given in
(\ref{E:second})
has annihilator equal to $(\lambda_{1}\cap\lambda_{3}) +
(\lambda_{2}\cap\lambda_{3})$.
\end{thm}
\begin{proof} Let $W$ denote the annihilator of this form.
It is clear that $\lambda_{1}\cap\lambda_{3} \subset W $, also
$\lambda_{2}\cap\lambda_{3} \subset  W$. Hence
$(\lambda_{1}\cap\lambda_{3})$ + $(\lambda_{2}\cap\lambda_{3})\subset  W$.
 Now to prove the other containment, let $a\in W$, so $\langle a,b\rangle$
= 0 for all $b\in (\lambda_{1} + \lambda_{2}) \cap \lambda_{3}$. In other
words; if $a=a_{1}+a_{2}\in \lambda_{3}$ where $a_{1}\in\lambda_{1}$ and
$a_{2}\in\lambda_{2}$, then we have $\psi(a_{2},b)$ = 0. It follows that
$a_{2}\in \Ann((\lambda_{1} + \lambda_{2}$) $\cap \lambda_{3})$ in $V$. Using equations \eqref{ann} and \eqref{ann2},
we have that
\begin{align}
\notag \Ann((\lambda_{1} + \lambda_{2})\cap\lambda_{3}) 
\notag  & = \Ann (\lambda_{1} + \lambda_{2})+  \Ann(\lambda_{3})\\
\notag & = ( \Ann(\lambda_{1} ) \cap \Ann(\lambda_{2}) )+ \lambda_{3} \\
\notag  & =(\lambda_{1}\cap\lambda_{2})+ \lambda_{3} 
\end{align} 
Thus $a_{2}\in
(\lambda_{1}\cap\lambda_{2})+\lambda_{3}$. So we could write $a_{2}=c+d$
where $c\in\lambda_{1}\cap\lambda_{2}$ and $d\in\lambda_{3}$.  It follows
that $a=(a_{1}+c)+d$ where $a_{1}+c\in\lambda_{1} $ and $d\in \lambda_{3}$.
Now since we have $a_{2}, c\in\lambda_{2}$ we get $d\in\lambda_{2}$.  
Since $a, d\in\lambda_{3}$ we get $ a_{1}+c\in\lambda_{3}$. Hence
$d\in\lambda_{2}\cap\lambda_{3}$ and
$a_{1}+c\in\lambda_{1}\cap\lambda_{3}$.  So $a=(a_{1}+c)+d\in
(\lambda_{1}\cap\lambda_{3})+(\lambda_{2}\cap\lambda_{3})$. Thus
$W\subset (\lambda_{1}\cap\lambda_{3})+(\lambda_{2}\cap\lambda_{3})$. So
$W=(\lambda_{1}\cap\lambda_{3})+(\lambda_{2}\cap\lambda_{3})$, i.e. the
annihilator of the form $\langle\, , \rangle$ is equal to
$(\lambda_{1}\cap\lambda_{3})+(\lambda_{2}\cap\lambda_{3})$.
\end{proof}

\begin{prop} \label{R:firstr}
 For any pair of Lagrangian subspaces $\lambda_{1}$, and $\lambda_{2}$ we
have
\[ \dim (\lambda_{1}) = \dim (\lambda_{2})
\]
and
\begin{align}\label{E:newequation}
  \dim (\lambda_{1} + \lambda_{2}) \equiv \dim (\lambda_{1}\cap\lambda_{2})
\quad \bmod(2).
\end{align}
\end{prop}

\begin{proof}
 The first formula follows by reducing it to the nonsingular case and
\begin{align}
\notag
\dim (A) = \dim (V) - \dim (\Ann(A)).
\end{align}
We obtain the second congruence from
\begin{align} \label{E:third}
\dim (A + B) = \dim (A) + \dim (B) - \dim (A\cap B)
\end{align}
 and the first formula.
 \end{proof}
 
 \begin{cor} \label{cor}
\begin{align}\notag
 \mu(\lambda_{1},\lambda_{2},\lambda_{3})\equiv \dim ((\lambda_{1} +
\lambda_{2}) \cap
\lambda_{3})+\dim ((\lambda_{1}\cap\lambda_{3})+(\lambda_{2}\cap\lambda_{3}))
\quad \bmod(2).
\end{align}
\end{cor}
\begin{proof} Since the annihilator  of the form is   $(\lambda_{1}\cap\lambda_{3})$ + $(\lambda_{2}\cap\lambda_{3})$, 
it follows that the rank of the form is 
\[\dim ((\lambda_{1} + \lambda_{2}) \cap
\lambda_{3})-\dim ((\lambda_{1}\cap\lambda_{3})+(\lambda_{2}\cap\lambda_{3}))
.\]  The result follows from the fact that the signature and the rank of a nondegenerate form agree modulo two.
\end{proof}

\begin{proof}[ Proof of Theorem \ref{maincong}] By  equation (\ref{E:third}), we have
\begin{align} \notag
 \dim (\lambda_{1} + \lambda_{2} + \lambda_{3}) \equiv \dim (\lambda_{1}) +
\dim (\lambda_{2} + \lambda_{3}) + \dim (\lambda_{1} \cap (\lambda_{2} +
\lambda_{3})) \quad \bmod(2)
\end{align}
and also have
\begin{align} \notag
\dim (\lambda_{1}\cap\lambda_{2}\cap\lambda_{3})   \equiv  
\dim (\lambda_{1}\cap\lambda_{2}) + \dim (\lambda_{1}\cap\lambda_{3})  +  
 \dim ((\lambda_{1}\cap\lambda_{2}) + (\lambda_{1}\cap \lambda_{3}))
\quad \bmod(2).
\end{align}

 Hence by Theorem (\ref{T:firstt}), the left hand sides of these two congruences are congruent. So their right hand sides must be congruent as well:
 \begin{align}
\notag \dim (\lambda_{1}) + \dim (\lambda_{2} + \lambda_{3}) + \dim (\lambda_{1}
\cap( \lambda_{2} + \lambda_{3})) \equiv
\dim (\lambda_{1}\cap\lambda_{2}) + \dim (\lambda_{1}\cap\lambda_{3})\\
\notag + \dim ((\lambda_{1}\cap\lambda_{2}) + (\lambda_{1}\cap \lambda_{3}))
\quad \bmod(2).
\end{align}
 The last equation is equivalent to
\begin{align}
\notag
\dim (\lambda_{1} \cap( \lambda_{2} + \lambda_{3})) +
\dim ((\lambda_{1}\cap\lambda_{2}) + (\lambda_{1}\cap \lambda_{3})) \equiv
\dim (\lambda_{1}) + \dim (\lambda_{2} + \lambda_{3})\\ \notag +
\dim (\lambda_{1}\cap\lambda_{2}) + \dim (\lambda_{1}\cap\lambda_{3})   \quad
\bmod(2).
\end{align}
 The left hand side of this last equation  is congruent to the Maslov index by Corollary \ref{cor}, and 
hence the first formula follows. The second formula follows
by equation \eqref{E:newequation}.
\end{proof}

\section {the weighted cobordism categories}
 
 All 3-manifolds and surfaces in this paper are assumed to be oriented and
compact.  We define a weighted cobordism category $\fC$ whose
objects are surfaces $\Si$  without boundary equipped with a Lagrangian subspace
$\lambda  \subset H_{1}(\Sigma;\BR)$. We will denote objects by pairs $(\Si, \lambda)$.
 A cobordism  from $(\Sigma ,\lambda) $ to  $(\Sigma',\lambda')$ is  a  3-manifold together with an orientation preserving homeomorphism (called its boundary identification) from its boundary to   $-\Sigma\sqcup\Sigma'$. Here, and
elsewhere,  $-\Sigma$ denotes $\Sigma$ with the opposite orientation.
Two cobordisms are equivalent if there is an orientation preserving homeomorphism between the underlying 3-manifolds that commutes with the boundary identifications.  A morphism  $M: (\Sigma ,\lambda)  \ra (\Sigma',\lambda')$ is  an equivalence class of cobordisms  from $(\Sigma ,\lambda) $ to  $(\Sigma',\lambda')$ together with an integer weight.  We denote morphisms by a single letter. We let $w(M)$ denote the weight of $M$.  We let $\flat M$ denote the underlying 3-manifold of a representative cobordism. This is well defined up to homeomorphism respecting the boundary identifications.   We call $(\Sigma ,\lambda)$ the source of $M$ and $(\Sigma' ,\lambda' )$ the target of $M$. We let $j_{M}$ denote  
the inclusion $\Sigma$ into  $\flat M,$ and $j^{M}$ denote the inclusion $\Sigma'$ into $\flat M.$ 
Here and sometimes below we ignore the boundary identifications for simplicity and we write as if $\Sigma \coprod \Sigma'$ were the boundary of $\flat M$.

We let ${M}_{*}(\lambda)$  denote the Lagrangian subspace \cite[p188-189]{T} of  $H_1(\Sigma',\BR)$ given by $ \left({ j^M_*}\right) ^{-1} \left( {{j_{M}}_ *}(\lambda) \right)$. Similarly we have the Lagrangian subspace of $H_1(\Sigma,\BR)$:
$ M^{*}(\lambda')=  \left( {j_{M}}_*\right) ^{-1} \left( {j^M_ *}(\lambda') \right).$
 
  If $M_1: (\Sigma ,\lambda)  \ra (\Sigma',\lambda')$ and $M_2: (\Sigma' ,\lambda')  \ra (\Sigma',\lambda'')$ are two morphisms we define $\flat (M_2 \circ M_1)$ by gluing $\flat M_2$ to $\flat M_1$ by identifying the target of $M_1$ to the source of $M_2.$ The boundary of this new 3-manifold is equipped with a boundary identification in the obvious way. The weight of the composition is given by the formula \footnote{As in \cite{G}, we adopt the sign  convention of \cite{W} rather than \cite{T} for the sign of the Maslov index term in this formula. It makes no real difference for this paper.}.
 \begin{align}\label{wofcomp}
w(M_2 \circ M_1)= w(M_1) + w(M_2)
- \mu({M_{1}}_{*}(\lambda),\lambda',{{M}_{2}}^{*}(\lambda'')) 
\end{align}

 The identity $\id_{(\Sigma,\lambda)}: (\Sigma,\lambda) \ra  (\Sigma,\lambda)$ is given 
 by $ \Si  \times I$ with the weight zero and the standard boundary identification.  Any morphism  $C:(\Sigma,\lambda) \ra  (\Sigma,\lambda') $ with $ \Si  \times I$ as the  underlying 3-manifold, and with the standard  boundary identification will be called a pseudo-cylinder over $\Sigma$.

\begin{lem}\label{L:thirdl} Pseudo-cylinders are invertible in  $\fC.$
The inverse  of  $C:(\Sigma,\lambda) \ra  (\Sigma,\lambda') $ is  the 
pseudo-cylinder from $(\Sigma,\lambda')$ to  $(\Sigma,\lambda)$ with weight $-w(C)$. 
\end{lem}

\begin{proof} This follows immediately from the definitions. One needs that  the Maslov index vanishes when two of the three lagrangians coincides \cite[p183]{T}
\end{proof}

If we make the same definitions but using Lagrangians in $H_1(\Sigma, \BQ)$, we obtain the cobordism category $\cc$ studied in \cite{G}.  As $H_1(\Sigma, \BQ) \otimes \BR$ is naturally isomorphic to $H_1(\Sigma, \BR)$, a Lagrangian in  $H_1(\Sigma, \BQ)$ determines one in $H_1(\Sigma, \BR).$  A Lagrangian of $H_1(\Sigma, \BR)$ which arises in this way is called rational. In this way, we obtain  a functor $\cc \ra \fC.$  

\section{the even cobordism category}

We repeat a definition from \cite{G} except now we apply it to morphisms of 
$\fC$  instead of $\cc$. We denote $\beta_i(\flat M)$ by  $\beta_i(M).$

\begin{Def} \label{evendef}
A cobordism $M:(\Sigma,\lambda) \rightarrow (\Sigma',\lambda')$  of $\fC$ is even if and only if
\begin{align} \notag
w(M) \equiv \dim \left({j_{M}}_{*}(\lambda) + j^M_{*}(\lambda' ) \right)+\beta_1(M)+\beta_0(M)  +
\beta_0(\Sigma)+ \frac {\beta_1(\Sigma')} 2+ \epsilon(M) \quad \bmod(2)
\end{align}
where  $\epsilon (M)$ is one if exactly one of $\Sigma$ and $\Sigma
'$ is nonempty and otherwise $\epsilon (M)$ is zero.  If a cobordism is not even, it is called odd.
\end{Def}

We note that the inverse of an even pseudo-cylinder is even.

The first author  showed  that the composite of two even
morphisms of  $\cc$  is again even \cite[Theorem 7.2]{G}. The subcategory $\cc^+$ was defined to be the category with the same objects as $\cc$ but with only  even morphisms. In the rest of this section,
we generalize this result to morphisms in $\fC.$ Given this result,  we  define the subcategory $\fC^+$ to be the category with the same objects as $\fC$ but with only the even morphisms. We would also get a subcategory if we left the $\epsilon(N)$ term out of Definition \ref{evendef}. However the definition that we give is more natural from some points of view \cite{G}.

\begin{prop}
A pseudo-cylinder  $C:(\Sigma,\lambda) \ra  (\Sigma,\lambda') $ is  even if and only if  
\begin{align} \notag
w(C) \equiv \frac{\beta_{1}(\Sigma)}{2}+\dim (\lambda+\lambda')  \quad
\bmod(2)
\end{align}

\end{prop}
\begin{proof} Apply the definition above.
\end{proof}

\begin{lem}\label{L:secondl}
Let $M:(\Sigma,\lambda) \rightarrow (\Sigma',\lambda')$ be an even morphism. If
$C:(\Sigma, \hat \lambda) \rightarrow (\Sigma,  \lambda)$ 
 and $C':(\Sigma',\lambda') \rightarrow (\Sigma', \tilde \lambda)$ be an even pseudo-cylinders,
then $M \circ C$ and $C' \circ M$ are even. \end{lem}

\begin{proof}   We first show  that  $M \circ C$ is even. 
We need to show
\begin{align}\label{E:Mc}
w(M \circ C) \equiv
\dim (j_{M *}(\hat \lambda) + j_{*}^{M}(\lambda' ) )+\beta_1(M)+ \beta_0(M) + \beta_0(\Sigma)+ \frac {\beta_1(\Sigma_{1})} 2+ \epsilon(M)
 \bmod(2)
\end{align}

By Equation \ref{wofcomp},
\begin{align}\label{we}
w(M \circ C) \equiv w(M) + w(C) +\mu(\hat \lambda,\lambda,M^{*}(\lambda'))  \bmod(2)
\end{align}

By assumption, we have that:
\begin{align}\label{E:M}
w(M) \equiv \dim (j_{M *}(\lambda) + {j_{*}}^{M}(\lambda'))+\beta_1(M)+\beta_0(M)+ \beta_0(\Sigma)  + \frac{\beta_1(\Sigma_{1})} 2+
\epsilon(M)  \bmod(2)
\end{align}
and,
\begin{align}\label{E:C}
w (C) \equiv \frac{\beta_{1}(\Sigma)}{2}+\dim ( \hat \lambda+\lambda)  \quad
\bmod(2)
\end{align}
 So after we substitute  \eqref{we}, \eqref{E:M} and \eqref{E:C} into \eqref{E:Mc}, we
conclude that we need only prove:
\begin{align} \notag
\dim (j_{M *}(\hat \lambda) + {j_{*}}^{M}(\lambda' ) ) + 
\mu(\hat \lambda,\lambda, M^{*}(\lambda')) + \dim (\hat \lambda + \lambda) +
\frac{\beta_1(\Sigma)} 2 \equiv \\ \notag \dim (j_{M *}(\lambda) +
{j_{*}}^{M}(\lambda' ) ) \quad \bmod(2)\end{align}
 Given Theorem \ref{maincong},  this last congruence becomes:
 
\begin{align} \label{con2}
\dim (j_{M *}(\hat \lambda) + {j_{*}}^{M}(\lambda' ) ) + \dim (\hat \lambda + 
M^{*}(\lambda')) \equiv \dim (j_{M *}(\lambda) + { j_{*}}^{M}(\lambda' ) )\ + 
\\ \notag \dim (\lambda + M^{*}(\lambda')) \quad \bmod(2)
\end{align}
 
For any subspace $\delta$ of $H_1(M,\BR)$, we have that 
\[ 
\delta + M^*(\lambda')  = \left( j_{M *} \right) ^{-1} \left( j_{M *} ( \delta) + j_*^M (\lambda') \right)
\]
as
\[ 
 j_{M *} \left( \delta + M^*(\lambda') \right) =  j_{M *} ( \delta) + j_*^M (\lambda')
\]
and kernel of $j_{M*} $ is a subset of $M^*(\lambda').$
Thus we have that  $\dim (\delta + M^{*}(\lambda')) = \dim (j_{M*}(\delta) +
j^M_{*}(\lambda'))  + n$  where $n$ is the dimension of  kernel of $j_{M*} $.  Thus  both sides of  \eqref{con2}  are congruent to $n$. Hence, we obtain
(\ref{E:Mc}).

The proof that $C' \circ M$ is even follows formally from the first part, if we consider how the parity of a cobordism changes when we reverse the orientation of the underlying 3-manifold and reverse the roles of source and target.
\end{proof}

\begin{prop}\label{P:firstp}
If there are even pseudo-cylinders $C$ and $C'$ over $\Sigma$, and $\Sigma'$ such that $ C \circ M \circ C'$
is even, then
 $M$ is an even cobordism in $\fC$ from $(\Sigma,\lambda)$ to
$(\Sigma',\lambda')$.
\end{prop}

\begin{proof} 
 It follows by lemma (\ref{L:thirdl}) that we can factor $M$ as $C^{-1}\circ C \circ M \circ C' \circ C'^{-1}$. Hence $M$ is even by two applications of  Lemma (\ref{L:secondl}).
\end{proof}

\begin{thm}
The composition of two even morphisms of $\fC$ is again even.
\end{thm}

\begin{proof}  Let $M_{1}, M_{2}$ be two even
morphisms and adopt the notations associated to $M_1$ and $M_2$ in \S 3. 
We need to show that $M_{2}\circ M_{1}$ is an even cobordism.
It suffices to show  $C'' \circ M_{2}\circ M_{1}\circ
C$ is even for some even pseudo-cylinders over  $C$ and $C''$ over $\Sigma$ and
$\Sigma''$ with rational Lagrangians  for $\Sigma$ and $\Sigma''$ .
On the other hand we can write $M_{2}\circ M_{1}$ as  $ M_{2}\circ C' \circ {C'} ^{-1} \circ M_{1}$
where $C'$ is an even pseudo-cylinder over  $\Sigma'$ whose  the target  has a rational Lagrangian.
We have that 
\[C''\circ M_{2}\circ M_{1}\circ C= C'' \circ M_{2}\circ C' \circ {C'} ^{-1} \circ M_{1}\circ C =  N_2 \circ N_1\]
where $N_2=C'' \circ M_{2}\circ C'$ and $N_1={C'} ^{-1} \circ M_{1}\circ C.$
By  Lemma \ref{L:secondl},  $N_1$, $N_2$ are even morphisms.  By Theorem 7.2 in \cite{G}, it follows that   $N_2 \circ N_1$ is even. Hence $M_{2}\circ
M_{1}$ is even.
\end{proof}

\end{document}